\documentclass{amsart}
\usepackage{amsmath,amssymb}
\usepackage{amscd,latexsym}

\newcommand{\Kh}{{\mathcal K}}

\newcommand{\Mh}{{\mathcal M}}

\newcommand{\be}{\mathbf{1}}

\newcommand{\dr}{\mathrm{dr}\,}

\newcommand{\dimnuc}{\mathrm{dim_{nuc}}\,}

\newcounter{number}[section]

\newenvironment{nummer}{\refstepcounter{number}{\noindent\arabic{section}.\arabic{number}}}{}

\newcommand{\bn}{\noindent \begin{nummer} \rm}
\newcommand{\en}{\end{nummer}}

\newenvironment{thms}{\noindent {\sc Theorem:} \it}{}
\newenvironment{lms}{\noindent {\sc Lemma:} \it}{}
\newenvironment{props}{\noindent {\sc Proposition:} \it}{}
\newenvironment{dfs}{\noindent {\sc Definition:} \it}{}

\newenvironment{remss}{\noindent {\sc Remarks:}}{}

\newenvironment{nots}{\noindent {\sc Notation:} }{}

\newenvironment{nproof}{\noindent {\sc Proof:}}{\mbox{}\hfill 
\rule[-.2ex]{.25em}{1.8ex}}

\begin{document}

\title[Nuclear dimension and corona factorization]{Nuclear dimension\\
 and the \\
 corona factorization property}

\author{Ping Wong Ng}
\address{Mathematics Department\\
University of Louisiana at Lafayette\\
217 Maxim D.\ Doucet Hall\\
Lafayette, LA 70504-1010\\
 USA}

\email{png@louisiana.edu}

\author{Wilhelm Winter}
\address{School of Mathematical Sciences\\
University of Nottingham\\ University Park\\
Nottingham NG7 2RD\\
United Kingdom}

\email{wilhelm.winter@nottingham.ac.uk}

\date{\today}
\thanks{{\it Supported by:} Louisiana State Board of Regents Grant 
LEQSF(2007-10)-RD-A35 and EPSRC \indent First Grant EP/G014019/1}

\setcounter{section}{-1}

\begin{abstract}
We show that stabilizations of  sufficiently noncommutative separable unital $C^{*}$-algebras with finite nuclear dimension have the corona factorization property. 
\end{abstract}

\maketitle

\section{Introduction}

\noindent
The classification of nuclear $C^{*}$-algebras has seen rapid progress in recent years, especially since the systematic exploitation of algebraic, topological and homological regularity properties, such as $\mathcal{Z}$-stability (where $\mathcal{Z}$ denotes the Jiang--Su algebra), finite topological dimension (most notably finite decomposition rank and finite nuclear dimension) and unperforation of the Cuntz semigroup. We refer to \cite{EllToms:regularity}, \cite{Winter:localizingEC}, \cite{Winter:dr-Z-stable} and \cite{TomsWinter:minhom} for an overview and a number of sample results.  It has been conjectured that the above mentioned properties are all equivalent, at least for separable, simple, unital and nuclear $C^{*}$-algebras. 

In this paper we are particularly interested in nuclear dimension; this is a notion of topological dimension for nuclear $C^{*}$-algebras which was introduced in \cite{WinterZac:dimnuc} and which generalizes the earlier concept of  decomposition rank, cf.\ \cite{KirWinter:dr}.  Finite decomposition rank has proven to be extremely useful for the classification of stably finite, simple, unital $C^{*}$-algebras. Nuclear dimension accesses substantially larger classes of $C^{*}$-algebras; for example, Kirchberg--Phillips classification has turned out to cover precisely those simple and traceless $C^{*}$-algebras which satisfy the universal coefficient theorem and have finite nuclear dimension, cf.\ \cite{WinterZac:dimnuc}. The present article contributes to the study of nuclear dimension and its relations to other notions of regularity for (nuclear) $C^{*}$-algebras.

The corona factorization property is a regularity property of stable $C^{*}$-algebras; while it can not in general be expected to be equivalent to any of the aforementioned, it is closely related to and was motivated by   
extension theory problems related to classification theory --- specifically, the theory of absorbing extensions.

The latter was originally developed to  
provide nice characterizations of $KK$-theory, and to prove 
the stable uniqueness theorems of classification theory. 
(See \cite{EllKuc:BDFabsorption}, \cite{Arv:ext}, \cite{Kas:ellK}, \cite{Kas:KK}, \cite{Voi:Weyl-vn}, \cite{DadEil:class}, \cite{Lin:stableau}.)

A major starting point was Elliott's and Kucerovsky's algebraic 
characterization of nuclearly absorbing extensions (see \cite{EllKuc:BDFabsorption}).
Basically, a separable stable $C^*$-algebra has the corona factorization
property if and only if it has abundantly many absorbing extensions 
\cite{KucNg:au}.  Among other things, this includes  many of the 
extensions that have been important in classification theory, and 
leads to simple and nice characterizations of $KK$-theory \cite{KucNg:au}.  
More recently, the corona factorization property has been used to classify
many nonsimple nuclear $C^*$-algebras (including interesting
examples coming from dynamical systems and graph theory; see,
for example, \cite{EilResRui:extensions}).  

 It turns out that the corona factorization property is also intimately
connected with the structure theory of $C^*$-algebras \cite{KucNg:S}, 
\cite{Ng:cfpsurvey}.   
Fundamental questions, like whether the extension of a separable stable
$C^*$-algebra by a separable stable $C^*$-algebra has a stable extension
algebra, are closely related to the corona factorization property
and the techniques used to study it.  

In the present article we study stabilizations of separable unital $C^{*}$-algebras with finite nuclear dimension; our main result (cf.\ Theorem~\ref{main-result-thm}) says that such algebras always have the corona factorization property, at least if they are sufficiently noncommutative:

\bn
\label{main-intro}
\begin{thms}
Let $A$ be a separable and unital $C^{*}$-algebra no hereditary subalgebra of which has a nontrivial elementary quotient. Suppose   $\dimnuc A \le n <\infty$. Then, $A \otimes \mathcal{K}$ has the corona factorization property.
\end{thms} 
\en

We remark that Ortega and R{\o}rdam have, independently from us and with different methods, arrived at a result similar to ours in \cite{OrtRor:cfp}: they show that stabilizations of separable, unital $C^{*}$-algebras with finite decomposition rank have the corona factorization property. Our Theorem~\ref{main-intro} partially generalizes this statement, since decomposition rank dominates nuclear dimension, cf.\ \cite{WinterZac:dimnuc}. The proof in \cite{OrtRor:cfp} heavily relies on systematic use of the Cuntz semigroup, and on $n$-comparison for $C^{*}$-algebras with decomposition rank at most $n$ as derived in \cite[Lemma~6.1]{TomsWinter:VI}. The respective statement is not known for $C^{*}$-algebras with finite nuclear dimension, which is the reason why we employ Kirchberg's covering number, cf.\ Section~{\ref{dimnuc}}. 

One of our motivations is to make progress on the question whether finite nuclear dimension implies $\mathcal{Z}$-stability for simple $C^{*}$-algebras. The respective statement for decomposition rank has been shown in \cite{Winter:dr-Z-stable}; that proof  also relied on \cite[Lemma~6.1]{TomsWinter:VI}. We are optimistic that our use of Kirchberg's covering number in place of $n$-comparison will contribute to the solution of the above mentioned question.

It is interesting to look at Theorem~\ref{main-intro} in the special case of simple $C^{*}$-algebras.
In this case, the result implies that stabilizations of all separable, simple, unital, nuclear $C^{*}$-algebras that have been classified so far do have the corona factorization property. This fact is not new, since all such $C^{*}$-algebras are known to be $\mathcal{Z}$-stable; however, our argument does not distinguish between the stably finite and the purely infinite case, hence points towards the possibility of a  unified approach to the classification problem for purely finite and purely infinite nuclear $C^{*}$-algebras.

The paper is organized as follows. We fix some notation in Section~{\ref{preliminaries}}. In Section~{\ref{cfp}} we recall the corona factorization property and some related facts on stable $C^{*}$-algebras. Section~{\ref{dimnuc}} states the definitions of nuclear dimension and of Kirchberg's covering number, and relates  the two concepts. Our main result is stated and proven in Section~{\ref{main-result}}.

\section{Preliminaries}
\label{preliminaries}

\bn
\label{def-subequivalence}
\begin{nots}
Let $A$ be a $C^{*}$-algebra and let $a,b \in A_{+}$ be positive elements. We write 
\[
a \preceq b, 
\]
if there is $x \in A$ such that 
\[
a= x^{*}x \mbox{ and } xx^{*} \le b. 
\]
We write $a \sim b$ if there is $x \in A$ such that
\[
a= x^{*}x \mbox{ and } b = xx^{*}.
\]
\end{nots}
\en

\bn
\label{rem-subequivalence}
\begin{remss}
(i) Note that, if in the situation of \ref{def-subequivalence} $a$ and $b$ are projections, then $a \preceq b$ means that $a$ is Murray--von Neumann subequivalent to $b$.

(ii) If $a$ is a projection and $b$ is a positive element of norm at most one, then $a \preceq b$ implies that for any $x \in A$ with $a=x^{*}x$ and $xx^{*} \le b$ we also have $bxx^{*} = xx^{*}$, hence $a = x^{*}bx$.
\end{remss}
\en

\bn
\begin{nots}
As usual, we denote by $M_{k}$ the $C^{*}$-algebra of complex $k$ by $k$ matrices, and by $\mathcal{K}$ the $C^{*}$-algebra of compact operators on a countably infinite dimensional Hilbert space. We will denote sets of matrix units for these algebras by $\{e_{ij} \mid i,j=1, \ldots,k\}$ and $\{f_{ij} \mid i,j=1,2,\ldots\}$, respectively.
\end{nots}
\en

\section{The corona factorization property}
\label{cfp}

\noindent
Below we recall the definition of the corona factorization property (cf.\ \cite{Ng:cfpsurvey}) and some related facts about stable $C^{*}$-algebras. We also derive a useful criterion for when a unital $C^{*}$-algebra has the corona factorization property.

\bn
\label{def-cfp}
\begin{dfs}
Let $B$ be a separable and stable $C^{*}$-algebra. $B$ is said to have the corona factorization property, if every full projection in $\mathcal{M}(B)$ is Murray--von Neumann equivalent to $\be_{\mathcal{M}(B)}$.
\end{dfs}
\en

\bn
\label{lem:full}
\begin{lms} 
 Let $A$ be a unital separable $C^*$-algebra.
Let $c \in \Mh(A \otimes \Kh)$ be a full element.
Then for every $b \in A \otimes \Kh$, $c + b$ is a full element of
$\Mh(A \otimes \Kh)$.
\end{lms}

\begin{nproof}
Since $c \in \Mh(A \otimes \Kh)$ is full, we may choose $x_1, x_2, .., x_n, 
y_1, y_2, ...,y_n \in \Mh(A \otimes \Kh)$ such that 
\begin{equation}
\label{ng1} \sum_{j=1}^{n} x_j c y_j = \be_{\Mh(A \otimes \Kh)}.
\end{equation}
Note that we have 
\begin{equation}
\label{ng2}
\sum_{i=1}^{\infty} \be_{A} \otimes f_{ii} = \be_{\Mh(A \otimes \Kh)},
\end{equation}
where the
sum converges in the strict topology on $\Mh(A \otimes \Kh)$. 

   Now since $A \otimes \Kh$ is an ideal of $\Mh(A \otimes \Kh)$, and
since $b \in A \otimes \Kh$, we have
$$b' := \sum_{i=1}^{n} x_i b y_i \in A \otimes \Kh.$$
Hence, by \eqref{ng2}, we can choose
$N \geq 1$ such that for all $m \geq N$,
\begin{equation}
\label{ng3}
\textstyle \left\| b' \left(\sum_{i \geq m} \be_{A} \otimes f_{ii}\right) \right\|, \; \left\| \left(\sum_{i \geq m} \be_{A} \otimes f_{ii} \right) b' \right\| < \varepsilon.
\end{equation}
Now let $s \in \Mh(A \otimes \Kh)$ be an isometry with range projection
$\sum_{i \geq N} \be_{A} \otimes f_{ii}$. 
Hence, by \eqref{ng3}, 
$$\| s^* b' s \| < \varepsilon.$$
Therefore, by \eqref{ng1},
\begin{eqnarray*}
\lefteqn{ \textstyle \left\| \sum_{j=1}^{n} s^* x_j (c + b) y_j s - \be_{\Mh(A \otimes \Kh)} \right\| }\\
& \leq & \textstyle \left\| \sum_{j=1}^{n} s^* x_j c y_j s - \be_{\Mh(A \otimes \Kh)} \right\| 
+ \| s^* b' s \| \\
& < & 0 + \varepsilon.
\end{eqnarray*}
Hence, if $\varepsilon > 0$ was chosen small enough, we see that 
$ \sum_{j=1}^{n} s^* x_j (c + b) y_j s$ must be invertible; and hence, 
$c + b$ is full in $\Mh(A \otimes \Kh)$.   
\end{nproof}
\en 

\bn
\label{lem:HjelmRor}
For a $C^*$-algebra $B$, let $F(B)$ denote the set of positive elements
which have a multiplicative identity, i.e.,
$$F(B) := \{ a \in B_+ \mid \exists b \in B_+ : ab = a \}.$$ 
We shall need the following result due to Hjelmborg and Rordam 
\cite{HjeRor:stable}:

\begin{lms} 
Let $B$ be a $\sigma$-unital $C^*$-algebra.
Then the following are equivalent:
\begin{enumerate}
\item $B$ is stable.
\item  For every $a \in F(B)$, there exists $b \in B_+$ such that 
$a \sim b$ and $a \perp b$. 
\end{enumerate}
\end{lms}
\en

\bn
\label{CFP_Aue}
The next result is \cite[Proposition~3.3]{KucNg:au}.

\begin{props} Let $B$ be a separable, stable $C^*$-algebra.   
Let $l$ be a nonzero positive element of $\mathcal{M}(B)$.  Then 
the hereditary subalgebra $\overline{l B l} \subset B$ is isomorphic to 
a hereditary subalgebra of the form $p B p$, where $p$ is a multiplier
projection.

   Moreover, if $l$ is a norm-full element of $\mathcal{M}(B)$, then $p$ can
also be chosen to be a norm-full element of $\mathcal{M}(B)$.  
\end{props}
\en

\bn
\label{B}
The following will be useful for the proof of our main result in Section~{\ref{main-result}}.

\begin{props}
Let $A$ be a separable and unital $C^{*}$-algebra. Suppose that, for any full projection $p \in \mathcal{M}(A \otimes \mathcal{K})$, we have $\be_{A} \otimes f_{11} \preceq p$.

Then, $A \otimes \mathcal{K}$ has the corona factorization property.
\end{props}

\begin{nproof}
Say that $p \in \Mh(A \otimes \Kh)$ is a full projection.
We will first show that $p(A \otimes \Kh)p$ is a stable $C^*$-algebra. 

   Let $a \in F(p(A \otimes \Kh)p)$ be given.  Hence, let $b \in 
p(A \otimes \Kh)p$ be a positive element such that $ab = a$.
For simplicity, let us assume that $\| a \| = \| b \| = 1$.

    By Lemma \ref{lem:full}, $p - b$ is a full positive element
of $\mathcal{M}(A \otimes \Kh)$.  Hence, by Proposition \ref{CFP_Aue}, 
let $q \in \mathcal{M}(A \otimes \Kh)$ be a full projection such that 
there exists a $*$-isomorphism
$$\phi : q(A \otimes \Kh)q \rightarrow \overline{(p - b) (A \otimes \Kh) 
(p- b)}.$$

    By hypothesis, let $f \in q(A \otimes \Kh)q$ be a projection such that
$f \sim \be_A \otimes f_{1,1}$.
Note that since $\be_A \otimes f_{1,1}$ is full in $A \otimes \Kh$,
$f$ is full in $q(A \otimes \Kh)q$.   
Hence, $\phi(f)$ is full in $\phi(q(A \otimes \Kh)q) = 
\overline{(p - b) (A \otimes \Kh) (p- b)}$. 
But since $p -b$ is full in $\mathcal{M}(A \otimes \Kh)$,
$\overline{(p - b) (A \otimes \Kh) (p- b)}$ is a full hereditary
subalgebra of $A \otimes \Kh$.      
Hence, $\phi(f)$ is a full element of 
$p (A \otimes \Kh)p$.  
Hence, there exists $N \geq 1$ such that 
$$b \preceq \bigoplus^N \phi(f).$$

   By repeated applications of Lemma \ref{lem:full}, by the hypothesis, 
and using a short induction
argument, we can show that      
there exists a projection $g \in q(A \otimes \Kh)q$ such that 
$$g \sim \bigoplus^N (\be_A \otimes f_{1,1}).$$

   Hence,  we have that 
$$a \preceq b \preceq \phi(g).$$
Moreover, since  
$\phi(g) \in \overline{(p - b) (A \otimes \Kh) (p- b)}$,
we also have that 
$$\phi(g) \perp a.$$ 

Since $a$ was arbitrary, it follows, by Lemma \ref{lem:HjelmRor}, that
$p (A \otimes \Kh)p$ is a stable $C^*$-algebra.  

From Brown's theorem \cite[Corollary~2.8]{Bro:stabher} we now see that $p(A \otimes \Kh)p \cong A \otimes \Kh$.
Hence, $$p \Mh(A \otimes \Kh) p \cong \Mh(p(A \otimes \Kh)p) \cong
\Mh(A \otimes \Kh).$$
(For the first isomorphism, see for example \cite[Lemma~11]{Kuc:ext}.) It follows that the unit of
$p \Mh(A \otimes \Kh) p$,
which is $p$, must be properly infinite. But it is straightforward to check that a properly infinite full projection is equivalent to the unit, so, $p \sim \be_{\Mh(A \otimes \Kh)}$ in $\Mh(A \otimes \Kh)$, as desired.
\end{nproof}
\en

\section{Nuclear dimension and Kirchberg's covering number}
\label{dimnuc}

\noindent
In this section we recall some concepts and results related to noncommutative topological dimension.

\bn
\label{def-dimnuc}
Recall the following definition from \cite{WinterZac:dimnuc}.

\begin{dfs}
A  $C^*$-algebra $A$ has nuclear dimension at most $n$, $\dimnuc A \le n$, if there exists a net 
$(F_{\lambda},\psi_{\lambda},\varphi_{\lambda})_{\lambda \in \Lambda}$ such that the $F_{\lambda}$ are finite-dimensional $C^{*}$-algebras, and such that $\psi_{\lambda}: A \to F_{\lambda}$ and  
$\varphi_{\lambda}: F_{\lambda} \to A$ are completely positive maps satisfying
\begin{enumerate}
\item $\varphi_{\lambda} \circ \psi_{\lambda} (a) \to a$ for any $a \in A$;
\item $\| \psi_{\lambda} \| \leq 1$; 
\item for each $\lambda$,  $F_{\lambda}$ decomposes into $n+1$ ideals $F_{\lambda}= F_{\lambda}^{(0)} \oplus \ldots \oplus F_{\lambda}^{(n)}$ 
such that $\varphi_{\lambda}|_{F_{\lambda}^{(i)}} $ is a c.p.c.\ order zero map for $i=0,1, \ldots ,n$.
\end{enumerate}
\end{dfs}
\en

\bn
\label{def-cov}
In \cite{Kir:CentralSequences}, Kirchberg introduced his notion of a covering number for a unital $C^{*}$-algebra. 

\begin{dfs}
Let $A$ be a unital $C^{*}$-algebra and $n \in \mathbb{N}$. $A$ has covering number at most $n$, $\mathrm{cov} A \le n$, if the following holds: 

For any $k \in \mathbb{N}$, there are a finite-dimensional $C^{*}$-algebra $F$, $d^{(1)}, \ldots,d^{(n)} \in A$  and a c.p.\ map $\varphi:F \to A$ such that
\begin{enumerate}
\item $F$ has no irreducible representation of rank less than $k$
\item $\varphi$ is $(n-1)$-decomposable with respect to $F = F^{(1)} \oplus \ldots \oplus F^{(n)}$
\item $\be_{A} = \sum_{j=1}^{n} (d^{(j)})^{*}\varphi^{(j)}(\be_{F^{(j)}})d^{(j)}$.
\end{enumerate} 
\end{dfs}
\en

\bn
\label{dimnuc-cov}
The following was shown in \cite[Proposition~4.3]{WinterZac:dimnuc}, the last statement  following from the proof of that result.

\begin{props}
Let $A$ be a separable $C^{*}$-algebra no  hereditary $C^{*}$-subalgebra of which has a nonzero elementary quotient, and let $\omega \in \beta \mathbb{N} \setminus \mathbb{N}$ be a free ultrafilter. Suppose $\dimnuc A \le n < \infty$. Then, 
\[
\mathrm{cov}(A_{\omega} \cap A'/\mathrm{Ann}(A)) \le (n+1)^{2},
\]
where
\[
\mathrm{Ann}(A):= \{b \in A_{\omega} \mid b A = A b = 0\}
\]
denotes the annihilator of $A$ in $A_{\omega}$.

In fact, for any $k \in \mathbb{N}$ and $l=1, \dots,(n+1)^{2}$, there are c.p.c.\ order zero maps
\[
\varrho^{(l)}:M_{k} \oplus M_{k+1} \to A _{\omega} \cap A '/\mathrm{Ann} (A) 
\]
such that
\[
\sum_{l=1}^{(n+1)^{2}} \varrho^{(l)}((\be_{k},\be_{k+1})) \ge \be_{A_{\omega} \cap A'/\mathrm{Ann}(A)}.
\]
\end{props}
\en

\bn
\label{C}
\begin{props}
Let $A$ be a $C^{*}$-algebra, $q \in A$ and $p \in \mathcal{M}(A)$ projections, and $\omega \in \beta \mathbb{N} \setminus \mathbb{N}$ a free ultrafilter. 

If $q \preceq \iota(p)$ in $\mathcal{M}(A)_{\omega}/ \mathrm{Ann}(A)$ (where $\iota:\mathcal{M}(A) \to  \mathcal{M}(A)_{\omega}/ \mathrm{Ann}(A)$ is the canonical embedding), then $q \preceq p$ in $\mathcal{M}(A)$. 
\end{props}

\begin{nproof}
We may lift the Murray--von Neumann relation $q \preceq \iota(p)$ to find a sequence of contractions 
\[
v = (v_{n})_{\mathbb{N}} \in \prod_{\mathbb{N}} \mathcal{M}(A)
\]
satisfying 
\[
q v_{n} p= v_{n} \mbox{ for all } n \in \mathbb{N} \mbox{ and } \lim_{\omega} (v_{n}v_{n}^{*} - q) \in \mathrm{Ann}(A).  
\]
Since $q \in A$, these relations in fact imply that
\[
\lim_{\omega} v_{n}v_{n}^{*} = q.
\]
Since $\omega$ is  free, there is a subsequence $(v_{n_{k}})_{k \in \mathbb{N}}$ such that 
\[
\lim_{k \to \infty} v_{n_{k}}v_{n_{k}}^{*} = q;
\]
note that we also have
\[
v_{n_{k}}^{*}v_{n_{k}} \le p \mbox{ for all } k \in \mathbb{N}.
\]
Now the assertion follows from the fact that Murray--von Neumann subequivalence is a semiprojective relation (cf.\ \cite{Lor:lifting}).
\end{nproof}
\en

\section{The main result}
\label{main-result}

\noindent
In this section, we derive our main result, Theorem~\ref{main-result-thm}.

\bn
\label{E}
\begin{props}
Let $A$ be a separable and unital $C^{*}$-algebra and $p \in \mathcal{M}(A \otimes \mathcal{K})$ a full positive contraction. Let $k \in \mathbb{N}$ and $q \in M_{k+1} \otimes A \otimes \mathcal{K}$ be a projection such that
\[
q \preceq \be_{k} \otimes p,
\]
where we have identified $M_{k} \otimes \mathcal{M} (A \otimes \mathcal{K})$ with its upper left corner embedding into $M_{k+1} \otimes \mathcal{M}(A \otimes \mathcal{K})$. 

Then, there are positive contractions
\[
d \in \overline{p (A \otimes \mathcal{K}) p} \mbox{ and } p' \in \overline{p \mathcal{M}(A \otimes \mathcal{K})p}
\]
such that $p'$ is full in $\mathcal{M}(A \otimes \mathcal{K})$, and such that
\[
d \perp p' \mbox{ and } q \preceq \be_{k} \otimes d.
\]
\end{props}

\begin{nproof}
Since $q$ is a projection in $M_{k+1} \otimes A \otimes \mathcal{K}$, it is straightforward to check that we may in fact assume that 
\[
q \in \be_{k} \otimes \overline{p (A \otimes \mathcal{K}) p} \subset M_{k} \otimes A \otimes \mathcal{K}.
\]
Let $(u_{n})_{\mathbb{N}}$ be an idempotent approximate unit for $\overline{p (A \otimes \mathcal{K})p}$, then
\[
(\be_{k} \otimes u_{n}^{\frac{1}{2}}) q (\be_{k} \otimes u_{n}^{\frac{1}{2}}) \to q.
\]
Again since $q$ is a projection, one checks that, if $n_{0}$ is large enough, then 
\[
q \preceq \be_{k} \otimes u_{n_{0}}.
\]
Since $(u_{n})_{\mathbb{N}}$ was idempotent, we have 
\[
u_{n_{0}} \perp (\be_{\mathcal{M}(A \otimes \mathcal{K}} - u_{n_{0}+1}) p (\be_{\mathcal{M}(A \otimes \mathcal{K})} - u_{n_{0}+1}).
\]
By Lemma~\ref{lem:full}, $(\be_{\mathcal{M}(A \otimes \mathcal{K})} - u_{n_{0}+1}) p (\be_{\mathcal{M}(A \otimes \mathcal{K})} - u_{n_{0}+1})$ is full in $\mathcal{M}(A \otimes \mathcal{K})$, so that
\[
d:= u_{n_{0}} 
\]
and 
\[
p' := (\be_{\mathcal{M}(A \otimes \mathcal{K})} - u_{n_{0}+1}) p (\be_{\mathcal{M}(A \otimes \mathcal{K})} - u_{n_{0}+1}) = p  (\be_{\mathcal{M}(A \otimes \mathcal{K})} - u_{n_{0}+1})^{2} p
\]
have the required properties.
\end{nproof}
\en

\bn
\label{A}
\begin{lms}
Let $A$ be a separable and unital $C^{*}$-algebra, $m \in \mathbb{N}$ and $p \in \mathcal{M}(A \otimes \mathcal{K})$ a full projection. 

Then, there are $ k \in \mathbb{N}$ and $d^{(1)}, \ldots, d^{(m)} \in (A \otimes \mathcal{K})_{+}$ with the following properties:
\begin{enumerate}
\item[(i)] the $d^{(l)}$ are pairwise orthogonal, have norm at most one and satisfy 
\[
\sum_{l=1}^{m} d^{(l)} \le p
\]
\item[(ii)] $\be_{k+1} \otimes \be_{A} \otimes f_{11} \preceq \be_{k} \otimes d^{(l)}$ for $l=1, \ldots,m$.
\end{enumerate}
\end{lms}

\begin{nproof}
Set
\begin{equation}
\label{ww1}
p^{(0)}:= p.
\end{equation}
Suppose that for some $0 \le l < m$, a full positive contraction $p^{(l)} \in \mathcal{M}(A \otimes \mathcal{K})$ has been constructed. Then, there is $k^{(l+1)} \in \mathbb{N}$ such that 
\[
e_{11} \otimes \be_{\mathcal{M}(A \otimes \mathcal{K})} \preceq \be_{k^{(l+1)}} \otimes p^{(l)}
\]
in $M_{k^{(l+1)}} \otimes \mathcal{M}(A \otimes \mathcal{K})$. But then also
\[
\be_{k^{(l+1)}+1} \otimes \be_{A} \otimes f_{11} \preceq \be_{k^{(l+1)}} \otimes p^{(l)}.
\]
Use Proposition~\ref{E} to find positive contractions
\[
d^{(l+1)} \in \overline{p^{(l)} (A \otimes \mathcal{K}) p^{(l)}} \mbox{ and } p^{(l+1)} \in \overline{p^{(l)} \mathcal{M}(A \otimes \mathcal{M}) p^{(l)}}
\]
such that $p^{(l+1)}$ is full in $\mathcal{M}(A \otimes \mathcal{K})$, and such that
\[
d^{(l+1)} \perp p^{(l+1)}
\]
and
\[
\be_{k^{(l+1)}+1} \otimes \be_{A} \otimes f_{11} \preceq \be_{k^{(l+1)}} \otimes d^{(l+1)}.
\]
Procced inductively to obtain $k^{(l)} \in \mathbb{N}$ and positive contractions $d^{(l)} \in A \otimes \mathcal{K}$ and $p^{(l)} \in \mathcal{M}(A \otimes \mathcal{K})$ satisfying
\begin{itemize}
\item[(a)] $d^{(l)} \in \overline{p^{(l-1)} (A \otimes \mathcal{K}) p^{(l-1)}}$
\item[(b)] $p^{(l)} \in \overline{p^{(l-1)} \mathcal{M}(A \otimes \mathcal{K}) p^{(l-1)}}$
\item[(c)] $p^{(l)}$ is full in $\mathcal{M}(A \otimes \mathcal{K})$
\item[(d)] $d^{(l)} \perp p^{(l)}$
\item[(e)] $\be_{k^{(l)}+1} \otimes \be_{A} \otimes f_{11} \preceq \be_{k^{(l)}} \otimes d^{(l)}$
\end{itemize}
for $l=1, \ldots,m$.

Note that (a), (b), (d) and \eqref{ww1} imply that the $d^{(l)}$ are positive contractions which are pairwise orthogonal and which satisfy
\[
\sum_{l=1}^{m} d^{(l)} \le p,
\]
so that \ref{A}(i) holds.

Set 
\[
k:= \prod_{i=1}^{m} k^{(i)} \mbox{ and } \check{k}^{(l)} := \frac{k}{k^{(l)}}
\]
for $l =1, \ldots,m$. We have (with the obvious identifications)
\begin{eqnarray*}
\be_{k+1} \otimes \be_{A} \otimes f_{11}& \le & \be_{k + \check{k}^{(l)}} \otimes \be_{A} \otimes f_{11} \\
& = & \be_{(k^{(l)}+1) \cdot \check{k}^{(l)}} \otimes \be_{A} \otimes f_{11} \\
& \stackrel{\mathrm{(e)}}{\le} & \be_{k^{(l)} \cdot \check{k}^{(l)}} \otimes d^{(l)}\\
& = & \be_{k} \otimes d^{(l)}
\end{eqnarray*}
for each $l \in 1, \ldots,m$, so \ref{A}(ii) holds.
\end{nproof}
\en

\bn
\label{main-result-thm}
We are now prepared to prove our main  result.

\begin{thms}
Let $A$ be a separable and unital $C^{*}$-algebra no hereditary $C^{*}$-subalgebra  of which has an elementary quotient. Suppose $\dimnuc A \le n <\infty$.

Then, $A \otimes \mathcal{K}$ has the corona factorization property.
\end{thms}

\begin{nproof}
Let $p \in \mathcal{M}(A \otimes \mathcal{K})$ be a full projection. In view of Proposition~\ref{B}, it will suffice to show that 
\begin{equation}
\label{w7}
\be_{A} \otimes f_{11} \preceq p \mbox{ in } \mathcal{M}(A \otimes \mathcal{K}).
\end{equation}
Set 
\[
m:= (n+1)^{2}
\]
and employ Lemma~\ref{A} to obtain $k \in \mathbb{N}$ and pairwise orthogonal elements 
\[
d^{(1)}, \ldots,d^{(m)} \in (A \otimes \mathcal{K})_{+}
\]
such that
\begin{equation}
\label{w20}
\sum_{l=1}^{m} d^{(l)} \le p
\end{equation}
and
\begin{equation}
\label{w8}
\be_{k+1} \otimes \be_{A} \otimes f_{11} \preceq \be_{k} \otimes d^{(l)} \mbox{ for } l=1, \ldots,m.
\end{equation}
Using \eqref{w8} and Remark~\ref{rem-subequivalence}(ii), we find 
\[
v_{k}^{(l)} \in M_{k} \otimes A \otimes \mathcal{K} \mbox{ and } v_{k+1}^{(l)} \in M_{k+1} \otimes A \otimes \mathcal{K}
\]
such that
\begin{equation}
\label{w13}
(v_{k}^{(l)})^{*}(\be_{k} \otimes d^{(l)}) v_{k+1}^{(l)} = \be_{k+1} \otimes \be_{A} \otimes f_{11}
\end{equation}
and
\[
(v_{k+1}^{(l)})^{*}(\be_{k} \otimes d^{(l)}) v_{k}^{(l)} = \be_{k} \otimes \be_{A} \otimes f_{11}
\]
Employ Proposition~\ref{dimnuc-cov} to find c.p.c.\ order zero maps
\begin{equation}
\label{w14}
\varrho^{(l)}: M_{k} \oplus M_{k+1} \to (A \otimes \mathcal{K})_{\omega} \cap (A \otimes \mathcal{K})'/\mathrm{Ann}(A \otimes \mathcal{K}) =: B
\end{equation}
for $l=1, \ldots,m$, such that 
\begin{equation}
\label{ww2}
\sum_{l=1}^{m} \varrho^{(l)}((\be_{k},\be_{k+1})) \ge \be_{B}.
\end{equation}
For $l=1, \ldots,m$, define linear maps
\[
\zeta_{k}^{(l)} : M_{k} \otimes M_{k} \to M_{k} \otimes B
\]
by
\begin{equation}
\label{w12}
e_{ij} \otimes e_{i'j'} \mapsto e_{ij} \otimes d^{(l)} \varrho^{(l)}((e_{i'j'},0)) \mbox{ for } i,j,i',j'=1, \ldots,k,
\end{equation}
and
\[
\zeta_{k+1}^{(l)} : M_{k+1} \otimes M_{k+1} \to M_{k+1} \otimes B
\]
by
\begin{equation}
\label{w16}
e_{ij} \otimes e_{i'j'} \mapsto e_{ij} \otimes d^{(l)} \varrho^{(l)}((0,e_{i'j'})) \mbox{ for } i,j,i',j'=1, \ldots,k+1.
\end{equation}
Using that
\begin{equation}
\label{w10}
[d^{(l)}, \varrho^{(l)}(M_{k} \oplus M_{k+1})] = 0,
\end{equation}
one checks that $\zeta_{k}^{(l)}$ and $\zeta_{k+1}^{(l)}$ are c.p.c.\ order zero maps and that 
\begin{equation}
\label{w9}
\zeta_{\bar{k}}^{(l)}(M_{\bar{k}} \otimes M_{\bar{k}}) \perp \zeta^{(l')}_{\tilde{k}}(M_{\tilde{k}} \otimes M_{\tilde{k}}) \mbox{ if } l \neq l' \in \{1, \ldots,m\} \mbox{ and } \bar{k},\tilde{k} \in \{k,k+1\}.
\end{equation}
Similarly, we may define c.p.c.\ order zero maps
\[
\xi_{k}^{(l)}: M_{k} \otimes M_{k} \to M_{k} \otimes B
\]
by
\begin{equation}
\label{w15}
e_{ij} \otimes e_{i'j'} \mapsto e_{ij} \otimes ((\be_{A} \otimes f_{11}) \varrho^{(l)}((e_{i'j'},0)))
\end{equation}
and
\[
\xi_{k+1}^{(l)}: M_{k+1} \otimes M_{k+1} \to M_{k+1} \otimes B
\]
by
\begin{equation}
\label{w15a}
e_{ij} \otimes e_{i'j'} \mapsto e_{ij} \otimes ((\be_{A} \otimes f_{11}) \varrho^{(l)}((0,e_{i'j'}))).
\end{equation}
Let
\[
\bar{\zeta}^{(l)}_{k}, \bar{\xi}^{(l)}_{k} :M_{k} \otimes M_{k} \to M_{k} \otimes B^{**}
\]
and
\[
\bar{\zeta}^{(l)}_{k+1}, \bar{\xi}^{(l)}_{k+1} :M_{k+1} \otimes M_{k+1} \to M_{k+1} \otimes B^{**}
\]
be the respective canonical supporting $*$-homomorphisms, cf.\ \cite{Winter:cpr2} and \cite{WinterZac:orderzero}. 

Next, choose partial isometries 
\[
s_{k} \in M_{k} \otimes M_{k} \mbox{ and } s_{k+1} \in M_{k+1} \otimes M_{k+1}
\]
satisfying
\begin{equation}
\label{w11}
s_{k} s_{k}^{*} = e_{11} \otimes \be_{k}, \, s_{k}^{*} s_{k} = \be_{k} \otimes e_{11}
\end{equation}
and
\begin{equation}
s_{k+1} s_{k+1}^{*} = e_{11} \otimes \be_{k+1}, \, s_{k+1}^{*} s_{k+1} = \be_{k+1} \otimes e_{11}.
\end{equation}
For $l=1, \ldots,m$, set
\begin{equation}
\label{w17}
u_{k}^{(l)} := \zeta^{(l)}_{k} (e_{11} \otimes \be_{k})^{\frac{1}{2}} \bar{\zeta}^{(l)}_{k}(s_{k}) v^{(l)}_{k} \bar{\xi}_{k}^{(l)}(s_{k}^{*}) \in M_{k} \otimes B^{**}
\end{equation}
and
\begin{equation}
\label{w18}
u_{k+1}^{(l)} := \zeta^{(l)}_{k+1} (e_{11} \otimes \be_{k+1})^{\frac{1}{2}} \bar{\zeta}^{(l)}_{k+1}(s_{k+1}) v^{(l)}_{k+1} \bar{\xi}_{k+1}^{(l)}(s_{k+1}^{*}) \in M_{k+1} \otimes B^{**}.
\end{equation}
We claim that in fact
\[
u_{k}^{(l)} \in M_{k} \otimes B \mbox{ and } u_{k+1}^{(l)} \in M_{k+1} \otimes B. 
\]
To this end, note that
\begin{eqnarray*}
u_{k}^{(l)} & = & \zeta_{k}^{(l)}(\be_{k} \otimes \be_{k})^{\frac{1}{2}} \bar{\zeta}^{(l)}_{k}((e_{11} \otimes \be_{k}) s_{k}) v_{k}^{(l)} \bar{\xi}^{(l)}_{k}(s_{k}^{*})  \\
& = & (\be_{k} \otimes d^{(l)})^{\frac{1}{4}} (\be_{k} \otimes d^{(l)} \varrho^{(l)}((\be_{k},0)))^{\frac{1}{4}}  (\be_{k} \otimes  \varrho^{(l)}((\be_{k},0)))^{\frac{1}{4}}  \\
&& \bar{\zeta}^{(l)}_{k} ((e_{11} \otimes \be_{k})s_{k}) v_{k}^{(l)} \bar{\xi}_{k}^{(l)}(s_{k}^{*})  \\
& = & (\be_{k} \otimes d^{(l)})^{\frac{1}{4}} (\zeta^{(l)}_{k})^{\frac{1}{4}} ((e_{11} \otimes \be_{k})s_{k}) v_{k}^{(l)} (\xi_{k}^{(l)})^{\frac{1}{4}}(s_{k}^{*}) \in M_{k} \otimes B, 
\end{eqnarray*}
where we have used functional calculus for order zero maps (cf.\ \cite{Winter:localizingEC} and \cite[3.2]{WinterZac:orderzero}) as well as the fact that the image of $\varrho^{(l)}$ commutes with $A \otimes \mathcal{K}$ in $B$. That $u_{k+1}^{(l)} \in M_{k+1} \otimes B$ is checked similarly. 

Next, we compute for $l,l' \in \{1, \ldots, m\}$
\begin{eqnarray}
(u_{k}^{(l')})^{*}u_{k}^{(l)} & \stackrel{\eqref{w17}}{=} & \bar{\xi}^{(l')}_{k}(s_{k}) (v_{k}^{(l')})^{*} \bar{\zeta}^{(l')}_{k}(s_{k}^{*}) \zeta^{(l')}_{k}(e_{11} \otimes \be_{k})^{\frac{1}{2}} \nonumber \\
&& \zeta^{(l)}_{k}(e_{11} \otimes \be_{k})^{\frac{1}{2}}  \bar{\zeta}^{(l)}_{k}(s_{k}) v_{k}^{(l)} \bar{\xi}^{(l)}_{k}(s_{k}^{*})  \nonumber \\
& \stackrel{\eqref{w9}}{=} & \delta_{ll'} \cdot \bar{\xi}_{k}^{(l)}(s_{k})(v_{k}^{(l)})^{*} \bar{\zeta}^{(l)}_{k}(s_{k}^{*}) \zeta^{(l)}_{k}(e_{11} \otimes \be_{k}) \bar{\zeta}_{k}^{(l)} (s_{k}) v_{k}^{(l)} \bar{\xi}_{k}^{(l)}(s_{k}^{*}) \nonumber \\
& = & \delta_{ll'} \cdot \bar{\xi}_{k}^{(l)}(s_{k})(v_{k}^{(l)})^{*}  \zeta^{(l)}_{k}(s_{k}^{*}(e_{11} \otimes \be_{k}) s_{k})  v_{k}^{(l)} \bar{\xi}_{k}^{(l)}(s_{k}^{*}) \nonumber \\
& \stackrel{\eqref{w11}}{=} &  \delta_{ll'} \cdot \bar{\xi}_{k}^{(l)}(s_{k})(v_{k}^{(l)})^{*}  \zeta^{(l)}_{k}(\be_{k} \otimes e_{11})  v_{k}^{(l)} \bar{\xi}_{k}^{(l)}(s_{k}^{*})  \nonumber \\
& \stackrel{\eqref{w12}}{=} &  \delta_{ll'} \cdot \bar{\xi}_{k}^{(l)}(s_{k})(v_{k}^{(l)})^{*}  (\be_{k} \otimes d^{(l)} \varrho^{(l)}((e_{11},0)))  v_{k}^{(l)} \bar{\xi}_{k}^{(l)}(s_{k}^{*})  \nonumber \\
& \stackrel{\eqref{w13},\eqref{w14}}{=} &  \delta_{ll'} \cdot \bar{\xi}_{k}^{(l)}(s_{k})  (\be_{k} \otimes ((\be_{A} \otimes f_{11})  \varrho^{(l)}((e_{11},0))))   \bar{\xi}_{k}^{(l)}(s_{k}^{*})  \nonumber \\
& \stackrel{\eqref{w15}}{=} &  \delta_{ll'} \cdot \bar{\xi}_{k}^{(l)}(s_{k})  \xi_{k}^{(l)}(\be_{k} \otimes e_{11})   \bar{\xi}_{k}^{(l)}(s_{k}^{*})  \nonumber \\
& = & \delta_{ll'} \cdot \xi_{k}^{(l)}(s_{k}(\be_{k} \otimes e_{11})s_{k}^{*}) \nonumber \\
& \stackrel{\eqref{w11}}{=} & \delta_{ll'} \cdot \xi_{k}^{(l)}(e_{11} \otimes \be_{k}) \nonumber \\
& \stackrel{\eqref{w15}}{=} & \delta_{ll'} \cdot e_{11} \otimes ((\be_{A} \otimes f_{11}) \varrho^{(l)}((\be_{k},0))). \label{w5}
\end{eqnarray}
Similarly, one checks
\begin{equation}
\label{w4}
(u_{k+1}^{(l)})^{*} u_{k+1}^{(l)} = \delta_{ll'} \cdot e_{11} \otimes ((\be_{A} \otimes f_{11}) \varrho^{(l)}((0,\be_{k+1}))).
\end{equation}
Let us identify $M_{k}$ with its upper left corner embedding into $M_{k+1}$, so that $u^{(l)}_{k}$ and $u_{k+1}^{(l)}$ all live in $M_{k+1} \otimes B$. Using \eqref{w12}, \eqref{w16}, \eqref{w10}, \eqref{w17}, \eqref{w18} and the fact that $\varrho^{(l)}$ has order zero, one checks that 
\begin{equation}
\label{w19}
(u_{k+1}^{(l')})^{*} u_{k}^{(l)} = 0  \mbox{ for all } l,l' \in \{1, \ldots,m\}. 
\end{equation}
We are now ready to define 
\begin{equation}
\label{w3}
v:= \sum_{l=1}^{m} (u^{(l)}_{k} + u^{(l)}_{k+1}) (e_{11} \otimes \sum_{l=1}^{m} \varrho^{(l)}((\be_{k},\be_{k+1})))^{-\frac{1}{2}} \in M_{k+1} \otimes B,
\end{equation}
where the inverse is taken in $e_{11} \otimes B$ (this is possible by \eqref{ww2}).
We compute
\begin{eqnarray}
v^{*}v & \stackrel{\eqref{w3}, \eqref{w5}, \eqref{w4},\eqref{w19}}{=} & (e_{11} \otimes \sum_{l=1}^{m} \varrho^{(l)}((\be_{k},\be_{k+1})))^{-\frac{1}{2}} \nonumber \\
&& (e_{11} \otimes ((\be_{A} \otimes f_{11}) \sum_{l=1}^{m} \varrho^{(l)}((\be_{k},\be_{k+1}))))  \nonumber \\
&& (e_{11} \otimes \sum_{l=1}^{m} \varrho^{(l)}((\be_{k},\be_{k+1})))^{-\frac{1}{2}} \nonumber \\
& \stackrel{\eqref{w14}}{=} & e_{11} \otimes \be_{A} \otimes f_{11}. \label{w1}
\end{eqnarray}
Furthermore, we note that
\[
u_{k}^{(l)} (u_{k}^{(l)})^{*} \stackrel{\eqref{w17}}{\le} \zeta_{k}^{(l)}(e_{11} \otimes \be_{k}) \stackrel{\eqref{w12}}{\le} e_{11} \otimes d^{(l)} \stackrel{\eqref{w20}}{\le} e_{11} \otimes p,
\]
so 
\[
(e_{11} \otimes p) u_{k}^{(l)} = u_{k}^{(l)},
\]
and, similarly,
\[
(e_{11} \otimes p) u_{k+1}^{(l)} = u_{k+1}^{(l)}.
\]
This implies
\[
(e_{11} \otimes p) vv^{*} = vv^{*}
\]
and, since $vv^{*}$ is a contraction, 
\begin{equation}
\label{w2}
vv^{*} \le e_{11} \otimes p.
\end{equation}
Note that \eqref{w1} and \eqref{w2} in particular imply that $v \in e_{11} \otimes B$. Identifying $e_{11} \otimes B$ with $B$, we see from \eqref{w1} and \eqref{w2} that
\[
\be_{A} \otimes f_{11} \preceq \iota(p)
\]
in $\mathcal{M}(A \otimes \mathcal{K})_{\omega}/\mathrm{Ann}(A \otimes \mathcal{K})$, where $\iota: \mathcal{M}(A \otimes \mathcal{K}) \to \mathcal{M}(A \otimes \mathcal{K})_{\omega}/\mathrm{Ann}(A \otimes \mathcal{K})$ denotes the canonical embedding via central sequences. By Proposition~\ref{C}, this implies \eqref{w7}.
\end{nproof}
\en

\bn
\label{typeI}
It follows directly from \cite[Theorem~2.9]{KucNg:typeI} that separable type I $C^{*}$-algebras with finite decomposition rank have the corona factorization property, cf.\ \cite[Theorem~3.1]{Ng:cfpsurvey}. We remark in closing that essentially the same proof yields the respective result for nuclear dimension in place of decomposition rank.

\begin{thms}
Let $A$ be a separable type I $C^{*}$-algebra. Suppose $\dimnuc A \le n <\infty$.

Then, $A \otimes \mathcal{K}$ has the corona factorization property.
\end{thms}

\begin{nproof}
The proof is essentially the same as that of \cite[Theorem~2.9]{KucNg:typeI}.  Starting with a composition series $(J_{\alpha})$ for $A$ such that each $J_{\alpha+1}/J_{\alpha}$ has continuous trace, there the permanence properties of the decomposition rank implied that the continuous trace algebras $J_{\alpha+1}/J_{\alpha}$ have decomposition rank at most $\dr A$. The same reasoning works for the nuclear dimension by \cite[Proposition~2.9]{WinterZac:dimnuc}. But since nuclear dimension and decomposition rank agree (with covering dimension of the spectrum) for continuous trace algebras (see \cite[Corollary~2.10]{WinterZac:dimnuc}), one can proceed just as in the proof of \cite[Theorem~2.9]{KucNg:typeI} from here. We omit the details. 
\end{nproof}
\en


\providecommand{\bysame}{\leavevmode\hbox to3em{\hrulefill}\thinspace}
\providecommand{\MR}{\relax\ifhmode\unskip\space\fi MR }
\providecommand{\MRhref}[2]{%
  \href{http://www.ams.org/mathscinet-getitem?mr=#1}{#2}
}
\providecommand{\href}[2]{#2}

\end{document}